\documentclass{article} 
\usepackage{mathmag-without-qedequation}

\usepackage{amsmath,amsthm}     
\usepackage{graphicx}     
\usepackage{hyperref} 
\usepackage{url}
\usepackage{amsfonts} 

\usepackage{multicol}

\usepackage[final]
{showkeys}

\usepackage[%
leftmargin=1cm,rightmargin=0\leftmargin
]
{quoting}

\theoremstyle{theorem}

\theoremstyle{definition}

\newtheorem{example}{Example}

\allowdisplaybreaks

\makeatletter
\@addtoreset{footnote}{page}
\makeatother

\renewcommand{\vec}{\mathbf}

\newcommand{\R}{\mathbb{R}}

\newcommand{\Ga}{\Gamma}
\newcommand{\de}{\delta}
\newcommand{\vp}{\varepsilon}

\newcommand{\cs}[1]{\langle #1 \rangle}

\newcommand{\curl}{\operatorname{curl}}

\begin{document}

\title{A simple and more general approach to Stokes' theorem}

\author{Iosif Pinelis\\               
\scriptsize Michigan Technological University\\    
ipinelis@mtu.edu}                      

\maketitle

\noindent  
{\small{\bf Summary } Oftentimes, Stokes' theorem is derived by using, more or less explicitly, the invariance of the curl of the vector field with respect to translations and rotations. However, this invariance -- which is oftentimes described as the curl being a ``physical'' vector -- does not seem quite easy to verify, especially for undergraduate students. 
An even bigger problem with Stokes' theorem is to rigorously define such notions as ``the boundary curve remains to the left of the surface''. 
Here an apparently simpler and more general approach is suggested.} 

\medskip

In 
calculus texts (see e.g.\ \cite[\S XII.6]{lang-UG} or \cite[\S 16.8]{stewart}), Stokes' theorem is usually stated as follows: Let $S$ be an oriented smooth enough surface in $\R^3$ bounded by a simple closed smooth enough curve $C$ with positive orientation. Then for any smooth enough vector field $\vec F$
\begin{equation}\label{eq:stokes-usual}
	\oint_C\vec F\cdot d\vec r=
	\iint_S\curl\vec F\cdot\vec n\,dS,  
\end{equation}
where $\curl\vec F$ is the curl of the vector field $\vec F$ and $\vec n$ is a continuous field of unit normal vectors on $S$. The positive orientation of $C$ is ``defined'' there as the condition for 
the surface $S$ to remain on the left (with respect to $\vec n$) when $C$ is traced out. 
In such elementary texts, it is of course impossible to rigorously define such notions as ``remains on the left''; cf.\ e.g.\  \cite[Theorems~5.5 and 5.9]{spivak}, \cite[\S XXIII.4--XXIII.6]{lang-RFA}, or \cite[\S 4.5.6]{federer}. The very concept of a surface in general is not elementary. A surface may be defined as a two-dimensional manifold-with-boundary \cite{spivak} -- possibly with singularities, which need to be considered to cover the case of even such simple surfaces as the (convex hulls of) triangles or rectangles; cf.\ e.g.\ \cite[\S XXIII]{lang-RFA}. Moreover, for the surface integral on the right-hand side of identity \eqref{eq:stokes-usual} to make sense, one has to ensure that the curl and the unit normal vector field have appropriate invariance properties with respect to the choice of an atlas for the manifold $S$.  

Green's theorem is the special and comparatively very simple case of Stokes' theorem corresponding to the additional condition that the surface $S$ is flat and thus may be assumed to coincide with a region $D\subset\R^2$. 
Informally, Green's theorem may be stated as follows: 

Let $D$ be a compact subset of $\R^2$ with a smooth enough boundary $\partial D$, and let $P$ and $Q$ be real continuously differentiable functions on an open neighborhood of $D$. Then 
\begin{equation}\label{eq:green}
	\oint_{\partial D}P\,du+Q\,dv=\iint_D \Big(\frac{\partial Q}{\partial u}-\frac{\partial P}{\partial v}\Big)du\, dv,  
\end{equation}
with the line integral $\oint_{\partial D}P\,du+Q\,dv$ appropriately defined. 

One way to make this statement rigorous (and general enough) is to assume that the boundary $\partial D$ of $D$ is a rectifiable Jordan curve $\Ga$ oriented so as to make the winding number of $\Ga$ equal $1$ (rather than $-1$) with respect to any point in the interior of $\Ga$. 
These conditions on $D$ and $\partial D$ will be assumed in the sequel. 

Even though Green's theorem is only a special case of Stokes', it is not easy to prove the just mentioned rigorous ``Jordan curve'' version of it, or even to show that the winding number can be consistently defined; see e.g.\ \cite{kostya} and references there to \cite{ridder,verblunsky,potts,apostol}. 
The best practical approach to teaching Green's theorem in a calculus course would probably be to restrict the consideration to simple regions, of the form $\{(x,y)\colon a\le x\le b, g_1(x)\le y\le g_2(x)\}$ or $\{(x,y)\colon a\le y\le b, g_1(y)\le x\le g_2(y)\}$, and possibly simple combinations thereof. 

Anyway, in this note we shall show that, once an appropriate version of Green's theorem is established, it is then very easy and painless to derive versions of Stokes' theorem, which are even more general, in some aspects, than the conventional one.  

Indeed, for some open neighborhood $U$ of $D$ 
and $(u,v) \in U$, let   
\begin{equation*}
(u,v)\longmapsto\vec r(u,v)
\in\R^3	
\end{equation*}
be a twice continuously differentiable map of 
$U$ into $\R^3$. 
Let $C:=\vec r(\partial D)$ be the image of the boundary $\partial D$ of $D$ under the map $\vec r$. 
Let 
\begin{equation*}
	\vec G\colon U\to\R^3
\end{equation*}
by a continuously differentiable vector field on $U$. Since $d\vec r=\vec r_u\,du+\vec r_v\,dv$, we can write/define the line integral $\oint_C\vec G\cdot d\vec r$ as follows: 
\begin{equation}\label{eq:oint_C}
	\oint_C\vec G\cdot d\vec r=\oint_{\partial D}(\vec G\cdot\vec r_u)\,du+(\vec G\cdot\vec r_v)\,dv. 
\end{equation}
Here, as usual, the subscripts ${}_u$, ${}_v$, ${}_x,\dots$ denote the partial differentiation with respect to $u,v,x,\dots$. 
One may say that formula \eqref{eq:oint_C} reduces the ``line integral'' $\oint_C\vec G\cdot d\vec r$ over the ``image-curve'' $C$ (which does not have to be a simple curve, without self-intersections) to the well-defined line integral $\oint_{\partial D}(\vec G\cdot\vec r_u)\,du+(\vec G\cdot\vec r_v)\,dv$ over the properly oriented rectifiable Jordan curve $\Ga=\partial D$ in the 
domain $U$ of the map $\vec r$.  

Letting now $P=\vec G\cdot\vec r_u$ and $Q=\vec G\cdot\vec r_v$, we have 
\begin{equation*}
	\frac{\partial Q}{\partial u}-\frac{\partial P}{\partial v}
	=\vec G_u\cdot\vec r_v+\vec G\cdot\vec r_{vu}-\vec G_v\cdot\vec r_u-\vec G\cdot\vec r_{uv}
	=\vec G_u\cdot\vec r_v-\vec G_v\cdot\vec r_u,  
\end{equation*}
since $\vec r_{vu}=\vec r_{uv}$. 
So, Green's theorem \eqref{eq:green} immediately yields the following form of Stokes' theorem: 
\begin{equation}\label{eq:stokes}
	\oint_C\vec G\cdot d\vec r=\iint_D \big(\vec G_u\cdot\vec r_v-\vec G_v\cdot\vec r_u\big)du\, dv,  
\end{equation} 
which, similarly to \eqref{eq:stokes-usual}, reduces a line integral to a double one.

This may be compared with the more conventional form of Stokes' theorem: 
\begin{equation}\label{eq:stokes-mine}
	\oint_C\vec F\cdot d\vec r=
	\iint_D(\vec\nabla\times\vec F)\cdot(\vec r_u\times\vec r_v)\,du\,dv,  
\end{equation}
where 
$\vec F$ is a continuously differentiable 3D vector field defined on a neighborhood of the ``surface'' $$S:=\vec r(D),$$
and the line integral may be understood as the one in \eqref{eq:stokes} (or, equivalently, in \eqref{eq:oint_C}) with 
\begin{equation}\label{eq:G=Fr}
	\vec G=\vec F\circ\vec r.   
\end{equation}


It is not hard to deduce \eqref{eq:stokes-mine} from \eqref{eq:stokes}. 
Indeed, let $(\vec i,\vec j,\vec k)$ be any orthonormal basis in $\R^3$, in which the cross products $\vec\nabla\times\vec F$ and $\vec r_u\times\vec r_v$ can be computed according to the standard determinant formulas. 
Let $\cs{x(u,v),y(u,v),z(u,v)}$ and $\cs{f(x,y,z),g(x,y,z),h(x,y,z)}$ be, respectively, the triples of the coordinates of  
$\vec r(u,v)$ and $\vec F(x,y,z)$ in the basis $(\vec i,\vec j,\vec k)$. 
Since both sides of \eqref{eq:stokes-mine} are linear in $\vec F$, without loss of generality $f=g=0$, and then the integrands on the right-hand sides of \eqref{eq:stokes} and \eqref{eq:stokes-mine} become, respectively,  
\begin{equation*}
\begin{aligned}
		I(u,v)&:=(h_x x_u+h_y y_u+h_z z_u) z_v-(h_x x_v+h_y y_v+h_z z_v) z_u \\ 
		&=(h_x x_u+h_y y_u)z_v-(h_x x_v+h_y y_v)z_u
\end{aligned}
\end{equation*}
and 
\begin{equation*}
\begin{aligned}
	J(u,v)&:=(h_y\vec i-h_x\vec j)\cdot\big((y_u z_v-z_u y_v)\vec i+(z_u x_v-x_u z_v)\vec j+(x_u y_v-y_u z_v)\vec  k\big) \\ 
	&=h_y(y_u z_v-z_u y_v)-h_x(z_u x_v-x_u z_v).  
\end{aligned}	
\end{equation*}
%
Now it is quite easy to see that $I(u,v)=J(u,v)$, which completes the derivation of \eqref{eq:stokes-mine} from \eqref{eq:stokes}. 

The main distinction of \eqref{eq:stokes} and \eqref{eq:stokes-mine} from \eqref{eq:stokes-usual} is that the double integrals in \eqref{eq:stokes} and \eqref{eq:stokes-mine} are taken over the flat preimage-region $D\subset\R^2$ -- whereas the ``image-surface'' $S=\vec r(D)$ plays no role in \eqref{eq:stokes} and \eqref{eq:stokes-mine}, except that the vector field $\vec F$ has to be defined on some neighborhood of $S$. 
%
Therefore, as far as \eqref{eq:stokes} and \eqref{eq:stokes-mine} are concerned, there is no need to talk about any properties of the image-surface $S$ except for it being a subset of $\R^3$. In particular, there is no need to talk about the orientation of $S$ or to use such hard to define terms as ``remains on the left'' or to care about the mentioned invariance properties of the curl and the unit normal vector field. 
Moreover, the image-surface $S=\vec r(D)$, as well as the image-curve $C=\vec r(\partial D)$, may be self-intersecting, and $S$ does not have to be a manifold at all. 
As for the line integrals in \eqref{eq:stokes} and \eqref{eq:stokes-mine}, they are over the image-curve $C$ only in form, as they immediately reduce to line integrals over the pre-image curve $\partial D$, according to \eqref{eq:oint_C}. 
This reduction of integration in the image-space $\R^3$ to that in the flat preimage-space $\R^2$ is  natural; it is even unavoidable -- for how else would one actually compute the line and surface integrals in the conventional form \eqref{eq:stokes-usual} of Stokes' theorem? 
(Of course, to do the integration in the preimage-space, 
we still need a proper orientation of the boundary $\partial D$ of the \emph{flat preimage} domain $D$, as was already pointed out in our discussion concerning a general ``Jordan curve'' version of Green's theorem and its elementary version for simple regions.)


Next, let us consider \eqref{eq:stokes} versus \eqref{eq:stokes-mine}. 
We saw that \eqref{eq:stokes} is very easy to obtain, modulo Green's theorem. Arguably, \eqref{eq:stokes} is also easier to remember than \eqref{eq:stokes-mine}. 

An important advantage of \eqref{eq:stokes} is that 
it is more general
than \eqref{eq:stokes-mine}. 
Indeed, on the one hand, 
\eqref{eq:stokes-mine} follows from \eqref{eq:stokes}; 
on the other hand, 
in \eqref{eq:stokes} the composition-factorization 
\eqref{eq:G=Fr} of the map $\vec G$ is not needed. That is, for \eqref{eq:stokes} one does not need the implication $\vec r(u_1,v_1)=\vec r(u_2,v_2)\implies \vec G(u_1,v_1)=\vec G(u_2,v_2)$ (which necessarily follows from \eqref{eq:G=Fr}). 

Plus, one does not need the notion of the curl for \eqref{eq:stokes}. On the other hand, \eqref{eq:stokes} by itself will not help when proving that a vector field is conservative if its curl is zero. Yet, as shown above, \eqref{eq:stokes-mine} is rather easy to get from \eqref{eq:stokes}. 


We thus have three versions of Stokes' theorem, corresponding to the three formulas: 
\begin{itemize}
	\item the conventional version \eqref{eq:stokes-usual}, which requires comparatively most stringent and even hard to define conditions, including an appropriate orientability of the image-surface together with the image-curve, and the invariance of the curl and the unit normal vector field; 
	\item the ``computational'', less conventional version \eqref{eq:stokes-mine}, which is not concerned with orientability of the images, but needs the composition-factorization condition \eqref{eq:G=Fr}; 
	\item \eqref{eq:stokes}, which is the most general of the three versions and, at the same time, easiest to obtain -- but not useful when, say, the curl is known to be zero. 
\end{itemize}

The above discussion 
is illustrated by

\begin{example}\label{ex:}
Consider the M\"obius strip $S_\de$ (of ``radius'' $1$ and half-width $\de>0$), which is the image $\vec r(D_\de)$ of the 
rectangle  
\begin{equation}\label{eq:D_de}
D_\de:=[0,2\pi]\times[-\de,\de] 	
\end{equation}
under the map $\vec r\colon D_\de\to\R^3$ given by the formula (cf.\ e.g.\ \cite{schwarz90})
\begin{equation}\label{eq:old}
\vec r(u,v)=\big((1+v\cos\tfrac u2)\cos u,(1+v\cos\tfrac u2)\sin u,v\sin\tfrac u2\big).  
\end{equation} 
%
The image-curve $C_\de:=\vec r(\partial D_\de)$ is self-intersecting, a reason being that for all $v\in[-\de,\de]$ one has $\vec r(0,v)=\vec r(2\pi,-v)[=(1+v,0,0)]$ -- whereas $(0,v)\in\partial D_\de$, $(2\pi,-v)\in\partial D_\de$, and $(0,v)\ne(2\pi,-v)$. 
For this reason, the ``surface'' $S_\de=\vec r(D_\de)$ may be considered self-intersecting as well. 
If $\de>2$, then $S_\de$ is self-intersecting in another, apparently more interesting manner; cf.\ \cite[Theorem~1.9]{moebius}. 
Indeed, assume that $\de>2$ and let $u_1=u$, $u_2=u+\pi$, $v_1=-2\cos\tfrac u2/\cos u$, and $v_2=-2\sin\tfrac u2/\cos u$, where $u$ is a small enough positive real number. Then 
$\vec r(u_1,v_1)=\vec r(u_2,v_2)[=(-1,-\tan u,-\tan u)\approx(-1,0,0)]$ -- whereas $(u_1,v_1)\in D_\de$, $(u_2,v_2)\in D_\de$, and $(u_1,v_1)\ne(u_2,v_2)$. 
Similarly, letting $u_1=2\pi-u$, $u_2=\pi-u$, $v_1=2\cos\tfrac u2/\cos u$, and $v_2=2\sin\tfrac u2/\cos u$  for small enough $u>0$, we have 
$\vec r(u_1,v_1)=\vec r(u_2,v_2)[=(-1,\tan u,\tan u)]$.

Letting now, for instance, $\vec G(u,v):=(u^2,0,0)$ for all $(u,v)$, we will have
$\vec G(u_1,v_1)\ne\vec G(u_2,v_2)$ for all pairs of points $(u_1,v_1)$ and $(u_2,v_2)$ as above. Therefore, at the self-intersection points $\vec r(u_1,v_1)=\vec r(u_2,v_2)=(-1,\vp\tan u,\vp\tan u)$ with $\vp=\pm1$ and small enough $u>0$, one will be unable to define a vector field $\vec F$ 
so that \eqref{eq:G=Fr} hold. 
Thus, formula \eqref{eq:stokes-mine} is not applicable here; of course, formula \eqref{eq:stokes-usual} is not applicable either, because the M\"obius strip is not orientable. 

In contrast, \eqref{eq:stokes} applies with no problem in this situation;   
%
each side of 
\eqref{eq:stokes} evaluates here to $-160\de/9$. 

The problem of application of Stokes' theorem to the M\"obius strip was previously considered in \cite{gabardo}. The version of the M\"obius strip dealt with in \cite{gabardo} differs by a composition of an isometry and a homothety from the apparently more common version described by \eqref{eq:old}. So, in the subsequent discussion the considerations in \cite{gabardo} will be translated into terms corresponding to \eqref{eq:old}. Following \cite{gabardo}, let us introduce here the vector field $\vec F$ defined by the formula 
\begin{equation}\label{eq:gab-F}
	\vec F(x,y,z)=\Big(\frac{-y}{x^2+y^2},\frac x{x^2+y^2},0\Big)
\end{equation}
for $(x,y,z)\in\R^3$ such that $x^2+y^2\ne0$. It is noted in \cite{gabardo} that (i) $\curl\vec F=\vec0$ wherever the vector field $\vec F$ is defined and (ii) the M\"obius strip does not intersect the $z$-axis if $\de<1$. So, $\curl\vec F=\vec0$ on the M\"obius strip. 

It is then concluded in \cite{gabardo} that the surface integral on the right-hand side of \eqref{eq:stokes-usual} is $0$. However, this conclusion is not quite correct, because the M\"obius strip is not orientable and hence the normal vector field $\vec n$ cannot be appropriately defined on the strip, so that the surface integral is technically not defined either and thus has no value. 


Also, it is observed in \cite{gabardo} that the boundary (say $B$) of the M\"obius strip is the image of the interval $[0,4\pi]$ under the map 
\begin{equation}\label{eq:tr}
u\mapsto\tilde{\vec r}(u):=\vec r(u,\de),  
\end{equation} 
and that the line integral $\oint_B\vec F\cdot d\tilde{\vec r}$ over the so-parameterized curve $B$ is $4\pi$, which differs from the presumed value $0$ of the actually undefined surface integral. 
This discrepancy is ascribed in \cite{gabardo} 
to the non-orientability of the M\"obius strip. 
Now one may wonder as follows: 

\begin{quoting}
We were told in this note that, as far as \eqref{eq:stokes} and \eqref{eq:stokes-mine} are concerned, there is no need to talk about the orientation of the images $S$ and $C$ of $D$ and $\partial D$ under the map $\vec r$. If so, then the version \eqref{eq:stokes-mine} of Stokes' theorem must hold even for the M\"obius strip and the vector field $\vec F$ as in \eqref{eq:gab-F}. The double integral on the right-hand side of \eqref{eq:stokes-mine} -- in contrast to that on the right-hand side of \eqref{eq:stokes-usual} -- is well defined, and it must be $0$, since $\curl\vec F=\vec0$ on the M\"obius strip. But the line integral over the boundary $B$ of the M\"obius strip was found in \cite{gabardo} to be nonzero. Does this not contradict \eqref{eq:stokes-mine}?  
\end{quoting}

In fact, there is no contradiction here. Recall that 
the line integral in \eqref{eq:stokes-mine} was to be understood as the one in \eqref{eq:oint_C} \big(with $\vec G$ as in \eqref{eq:G=Fr}\big) and, in turn, the line integral in \eqref{eq:oint_C} over the ``image-curve'' $C$ was defined as the corresponding line integral over the preimage $\partial D$ of $C$, where $\partial D$ is the boundary of the preimage-region $D$. 
In contrast with these conditions, the line integral in \cite{gabardo} was essentially taken over the segment $[0,4\pi]\times\{\de\}$, which 
is of course \emph{not} the boundary of $D=D_\de=[0,2\pi]\times[-\de,\de]$.  %
However, the value of the line integral in \eqref{eq:stokes-mine} computed indeed as the one in \eqref{eq:oint_C} with $\vec G$ as in \eqref{eq:G=Fr}
%
is $0$, which is of course the same as the value of the double integral in \eqref{eq:stokes-mine}.  

At this point, one may still wonder: 

\begin{quoting}
The difficulty with the applicability of Stokes' theorem to the M\"obius strip was ascribed in \cite{gabardo} to the non-orientability. 
However, the boundary $B$ of the (non-orientable) M\"obius strip
can also be the boundary of an orientable surface. Then the line integral $\oint_B\vec F\cdot d\tilde{\vec r}$ as computed in \cite{gabardo} will have the same value $4\pi\ne0$, whereas the corresponding surface integral on the right-hand side of \eqref{eq:stokes-usual} will still be $0$, right? So, it seems we have another contradiction here. 
\end{quoting}

The answer to this concern is as follows. 
First here, indeed the boundary $B$ of the M\"obius strip
can also be the boundary of an orientable surface. 
To see this, it is more convenient to use the interval $[-\pi,3\pi]$ instead of $[0,4\pi]$.  
More specifically, recall \eqref{eq:old} and note that, as $u$ increases from $-\pi$ to $\pi$, the $z$-coordinate $\de\sin\tfrac u2$ of the vector $\vec r(u,\de)$ increases from its minimal value $-\de$ to its maximal value $\de$; and as $u$ increases further from $\pi$ to $3\pi$, the $z$-coordinate of $\vec r(u,\de)$ decreases from $\de$ back to $-\de$. Moreover, 
each value of the $z$-coordinate of $\vec r(u,\de)$ in the interval $[-\de,\de)$ is taken at exactly two points $u_1(z)$ and $u_2(z)$ in $[-\pi,3\pi]$ such that $-\pi\le u_1(z)<\pi<u_2(z)=2\pi-u_1(z)\le3\pi$.  
The value $\de$ of the $z$-coordinate of $\vec r(u,\de)$ for $u\in[-\pi,3\pi]$ is taken only at $u=\pi$; accordingly, assume that $u_1(\de)=u_2(\de)=\pi$. 
Note also that $2\pi-u$ decreases from $3\pi$ to $\pi$ as $u$ increases from $-\pi$ to $\pi$. 
Connecting now, for each $z\in[-\de,\de]$, the points $\vec r(u_1(z),\de)$  and $\vec r(u_2(z),\de)$ by (say) a straight line segment, we do obtain an orientable surface, sat $\hat S$, whose boundary is the same as the boundary $B$ of the M\"obius strip. 
The surface $\hat S$ is the image $\hat{\vec r}(\hat D)$ of the rectangle 
$\hat D:=[-\pi,\pi]\times[0,1]$ under the map 
\begin{equation*}
	\hat D\ni(u,t)\mapsto\hat{\vec r}(u,t):=(1-t)\vec r(u,\de)+t\,\vec r(2\pi-u,\de).   
\end{equation*}
To verify that $\hat S$ is orientable, note that the first two coordinates 
of the vector $\vec c:=\hat{\vec r}_u(u,t)\times\hat{\vec r}_t(u,t)$ are $\de\cos\frac{u}{2}\,\sin u$ and $-\de^2 \cos^2\frac{u}{2}\, \cos u$, respectively, whence \break 
$\vec c\ne\vec0$ for all $(u,t)$ in the interior of $\hat D$, and so, one can define the unit normal vector field on 
the image of the interior of $\hat D$ under the map $\hat{\vec r}$ by the formula $\vec n=\vec c/|\vec c|$.  

Moreover and more importantly, the image $\hat{\vec r}(\partial\hat D)$ under this map of the boundary $\partial\hat D$ of the rectangle $\hat D$ is the boundary $B$ of the M\"obius strip. 
So, \eqref{eq:stokes-mine} will hold for any 3D vector field $\vec F$ which is smooth enough, in the sense of being continuously differentiable on a neighborhood of the surface $\hat S=\hat{\vec r}(\hat D)$. 
For instance, for the linear vector field defined by the formula $\vec F(x,y,z)=A(x\ y\ z)^\top$, where $A=(a_{i,j})_{i,j=1}^3$ is a constant $3\times3$ real matrix and ${}^\top$ denotes the transposition, both the left-hand side and right-hand side of \eqref{eq:stokes-mine} with $D=\hat D$ \big(and $C=\hat{\vec r}(\partial\hat D)$\big) take the same value, \break 
$2\pi(2+\de^2)(a_{1,2}-a_{2,1})+\pi\de^2(a_{1,3}-a_{3,1})$. 

Is the particular field $\vec F$ given by \eqref{eq:gab-F} continuously differentiable on a neighborhood of the surface $\hat S$? In other words, is the intersection of the surface $\hat S$ with the $z$-axis empty? 
If that were so, then the right-hand of \eqref{eq:stokes-mine} would be $0$ -- whereas, as one can recheck, the value of the line integral in \eqref{eq:stokes-mine} is the same nonzero value, $4\pi$, 
as the one found the other way in \cite{gabardo}. It follows that the surface $\hat S$ does intersect the $z$-axis. In fact, it is not hard to check directly that this intersection contains exactly two points, $(0,0,\pm\de/\sqrt2)$. 

A free bonus of this discussion is the fact that any orientable surface in $\R^3$ whose boundary coincides with the boundary of the M\"obius strip given by \eqref{eq:old} must necessarily intersect the $z$-axis. 
One may also note here that, for any such orientable surface and for $\vec F$ as in \eqref{eq:gab-F}, the surface integral on the right-hand side of \eqref{eq:stokes-usual} will actually not be $0$; rather, it will be undefined. 

This discussion is partly illustrated in Fig.~1. 

\begin{figure}[h]
	\centering
		\includegraphics[width=1.00\textwidth]
		{
		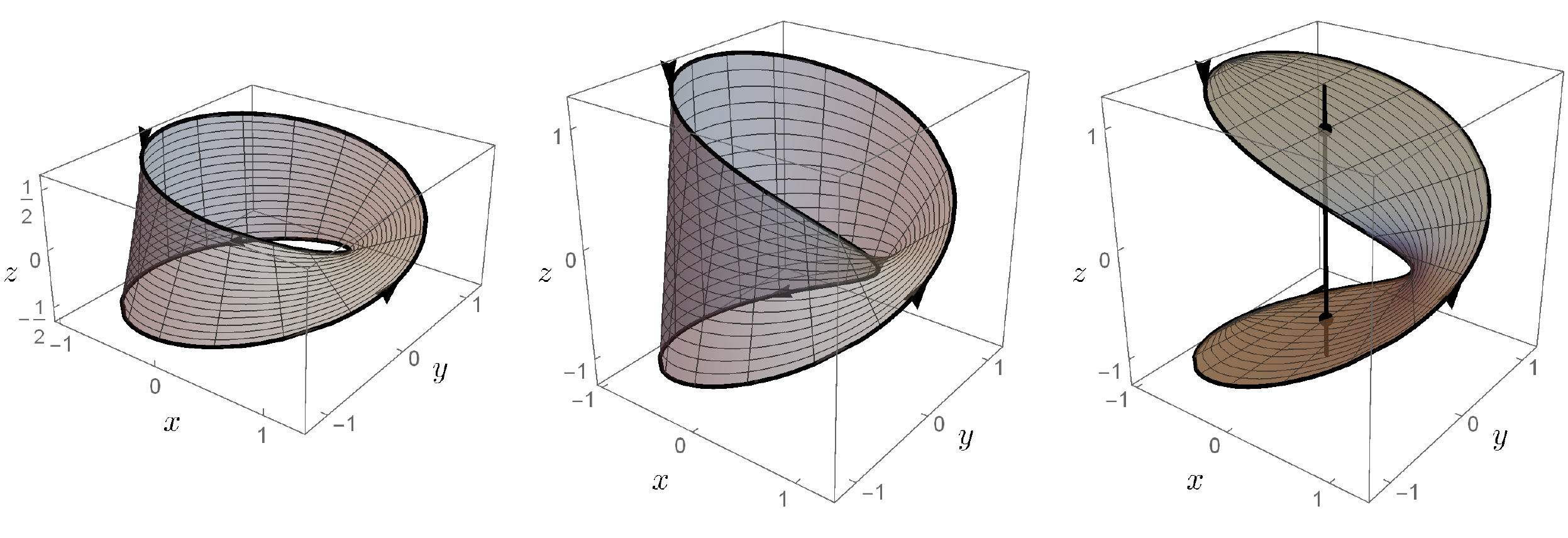
}
	\label{fig:}
	\caption{Left panel: the M\"obius strip given by \eqref{eq:old} (with $\de=3/10$), stretched vertically by the factor of $2$. Middle panel: the same M\"obius strip, now stretched vertically by the factor of $4$. Right panel: the oriented surface $\hat S$ with the same boundary as the M\"obius strip, stretched vertically by the factor of $4$, for better viewing; shown here are also the two points of intersection of the surface $\hat S$ with the $z$-axis.}  
\end{figure}

\end{example}

%

\def\cprime{$'$} \def\polhk#1{\setbox0=\hbox{#1}{\ooalign{\hidewidth
  \lower1.5ex\hbox{`}\hidewidth\crcr\unhbox0}}}
  \def\polhk#1{\setbox0=\hbox{#1}{\ooalign{\hidewidth
  \lower1.5ex\hbox{`}\hidewidth\crcr\unhbox0}}}
  \def\polhk#1{\setbox0=\hbox{#1}{\ooalign{\hidewidth
  \lower1.5ex\hbox{`}\hidewidth\crcr\unhbox0}}} \def\cprime{$'$}
  \def\polhk#1{\setbox0=\hbox{#1}{\ooalign{\hidewidth
  \lower1.5ex\hbox{`}\hidewidth\crcr\unhbox0}}} \def\cprime{$'$}
  \def\polhk#1{\setbox0=\hbox{#1}{\ooalign{\hidewidth
  \lower1.5ex\hbox{`}\hidewidth\crcr\unhbox0}}} \def\cprime{$'$}
  \def\cprime{$'$}

\noindent{\small {\bf IOSIF PINELIS }(MR Author ID: 208523) 
has been with Michigan Technological University since 1992. He has also held visiting positions at the University of Illinois, Urbana--Champaign; the City University of New York; and Lehigh University. 
His most extensive expertise is in probability and statistics, including
extremal problems, exact inequalities, and limit theorems of probability and statistics. He has also enjoyed doing work in machine learning, information theory, numerical analysis, operations research, combinatorics, mechanical engineering, biology, geometry, and physics.
}

%
%
%

\end{document}